\documentclass[12pt]{article}

\catcode`\@=11

\thicklines
\newskip\Einheit \Einheit=.6cm
\newcount\xcoord \newcount\ycoord
\newdimen\xdim \newdimen\ydim \newdimen\PfadD@cke \newdimen\Pfadd@cke
\def\PfadDicke#1{\PfadD@cke#1 \divide\PfadD@cke by2 
\Pfadd@cke\PfadD@cke \multiply\PfadD@cke by2}
\long\def\LOOP#1\REPEAT{\def\BODY{#1}\ITERATE}
\def\ITERATE{\BODY \let\next\ITERATE \else\let\next\relax\fi \next}
\let\REPEAT=\fi
\def\Punkt{\hbox{\raise-2pt\hbox to0pt{\hss\scriptsize$\bullet$\hss}}}

\def\DuennPunkt(#1,#2){\unskip
  \raise#2 \Einheit\hbox to0pt{\hskip#1 \Einheit
          \raise-1.5pt\hbox to0pt{\hss\tiny$\bullet$\hss}\hss}}
		  
\def\NormalPunkt(#1,#2){\unskip
  \raise#2 \Einheit\hbox to0pt{\hskip#1 \Einheit
          \raise-3pt\hbox to0pt{\hss\large$\bullet$\hss}\hss}}
\def\DickPunkt(#1,#2){\unskip
  \raise#2 \Einheit\hbox to0pt{\hskip#1 \Einheit
          \raise-4pt\hbox to0pt{\hss\Large$\bullet$\hss}\hss}}
\def\Kreis(#1,#2){\unskip
  \raise#2 \Einheit\hbox to0pt{\hskip#1 \Einheit
          \raise-4pt\hbox to0pt{\hss\Large$\circ$\hss}\hss}}
\def\Diagonale(#1,#2)#3{\unskip\leavevmode
  \xcoord#1\relax \ycoord#2\relax
      \raise\ycoord \Einheit\hbox to0pt{\hskip\xcoord \Einheit
         \unitlength\Einheit
         \line(1,1){#3}\hss}}
\def\AntiDiagonale(#1,#2)#3{\unskip\leavevmode
  \xcoord#1\relax \ycoord#2\relax \advance\xcoord by -0.05\relax
      \raise\ycoord \Einheit\hbox to0pt{\hskip\xcoord \Einheit
         \unitlength\Einheit
         \line(1,-1){#3}\hss}}
\def\Pfad(#1,#2),#3\endPfad{\unskip\leavevmode
  \xcoord#1 \ycoord#2 \thicklines\ZeichnePfad#3\endPfad\thinlines}
\def\ZeichnePfad#1{\ifx#1\endPfad\let\next\relax
  \else\let\next\ZeichnePfad
    \ifnum#1=1
      \raise\ycoord \Einheit\hbox to0pt{\hskip\xcoord \Einheit
         \vrule height\Pfadd@cke width1 \Einheit depth\Pfadd@cke\hss}%
      \advance\xcoord by 1
     \else\ifnum#1=2
      \raise\ycoord \Einheit\hbox to0pt{\hskip\xcoord \Einheit
         \unitlength\Einheit
         \line(0,1){1}\hss}
      \advance\xcoord by 0
      \advance\ycoord by 1
 \else\ifnum#1=3
      \raise\ycoord \Einheit\hbox to0pt{\hskip\xcoord \Einheit
         \unitlength\Einheit
         \line(1,1){1}\hss}
      \advance\xcoord by 1
      \advance\ycoord by 1
    \else\ifnum#1=4
      \raise\ycoord \Einheit\hbox to0pt{\hskip\xcoord \Einheit
         \unitlength\Einheit
         \line(1,-1){1}\hss}
      \advance\xcoord by 1
      \advance\ycoord by -1
   \else\ifnum#1=5
      \raise\ycoord \Einheit\hbox to0pt{\hskip\xcoord \Einheit
         \unitlength\Einheit
         \line(2,1){2}\hss}
      \advance\xcoord by 2
      \advance\ycoord by 1
	  \else\ifnum#1=6
      \raise\ycoord \Einheit\hbox to0pt{\hskip\xcoord \Einheit
         \unitlength\Einheit
         \line(2,-1){2}\hss}
      \advance\xcoord by 2
      \advance\ycoord by -1
	  \else\ifnum#1=7
      \raise\ycoord \Einheit\hbox to0pt{\hskip\xcoord \Einheit
         \unitlength\Einheit
         \line(3,1){3}\hss}
      \advance\xcoord by 3
      \advance\ycoord by 1
	  \else\ifnum#1=8
      \raise\ycoord \Einheit\hbox to0pt{\hskip\xcoord \Einheit
         \unitlength\Einheit
         \line(3,-1){3}\hss}
      \advance\xcoord by 3
      \advance\ycoord by -1
    \fi\fi\fi\fi\fi\fi\fi\fi
  \fi\next}
\def\hSSchritt{\leavevmode\raise-.4pt\hbox 
to0pt{\hss.\hss}\hskip.2\Einheit
  \raise-.4pt\hbox to0pt{\hss.\hss}\hskip.2\Einheit
  \raise-.4pt\hbox to0pt{\hss.\hss}\hskip.2\Einheit
  \raise-.4pt\hbox to0pt{\hss.\hss}\hskip.2\Einheit
  \raise-.4pt\hbox to0pt{\hss.\hss}\hskip.2\Einheit}
\def\vSSchritt{\vbox{\baselineskip.2\Einheit\lineskiplimit0pt
\hbox{.}\hbox{.}\hbox{.}\hbox{.}\hbox{.}}}
\def\DSSchritt{\leavevmode\raise-.4pt\hbox to0pt{%
  \hbox to0pt{\hss.\hss}\hskip.2\Einheit
  \raise.2\Einheit\hbox to0pt{\hss.\hss}\hskip.2\Einheit
  \raise.4\Einheit\hbox to0pt{\hss.\hss}\hskip.2\Einheit
  \raise.6\Einheit\hbox to0pt{\hss.\hss}\hskip.2\Einheit
  \raise.8\Einheit\hbox to0pt{\hss.\hss}\hss}}
\def\dSSchritt{\leavevmode\raise-.4pt\hbox to0pt{%
  \hbox to0pt{\hss.\hss}\hskip.2\Einheit
  \raise-.2\Einheit\hbox to0pt{\hss.\hss}\hskip.2\Einheit
  \raise-.4\Einheit\hbox to0pt{\hss.\hss}\hskip.2\Einheit
  \raise-.6\Einheit\hbox to0pt{\hss.\hss}\hskip.2\Einheit
  \raise-.8\Einheit\hbox to0pt{\hss.\hss}\hss}}
\def\SPfad(#1,#2),#3\endSPfad{\unskip\leavevmode
  \xcoord#1 \ycoord#2 \ZeichneSPfad#3\endSPfad}
\def\ZeichneSPfad#1{\ifx#1\endSPfad\let\next\relax
  \else\let\next\ZeichneSPfad
    \ifnum#1=1
      \raise\ycoord \Einheit\hbox to0pt{\hskip\xcoord \Einheit
         \hSSchritt\hss}%
      \advance\xcoord by 1
    \else\ifnum#1=2
      \raise\ycoord \Einheit\hbox to0pt{\hskip\xcoord \Einheit
        \hbox{\hskip-2pt \vSSchritt}\hss}%
      \advance\ycoord by 1
    \else\ifnum#1=3
      \raise\ycoord \Einheit\hbox to0pt{\hskip\xcoord \Einheit
         \DSSchritt\hss}
      \advance\xcoord by 1
      \advance\ycoord by 1
    \else\ifnum#1=4
      \raise\ycoord \Einheit\hbox to0pt{\hskip\xcoord \Einheit
         \dSSchritt\hss}
      \advance\xcoord by 1
      \advance\ycoord by -1
    \fi\fi\fi\fi
  \fi\next}
\def\Koordinatenachsen(#1,#2){\unskip
 \hbox to0pt{\hskip-.5pt\vrule height#2 \Einheit width.5pt depth1 
\Einheit}%
 \hbox to0pt{\hskip-1 \Einheit \xcoord#1 \advance\xcoord by1
    \vrule height0.25pt width\xcoord \Einheit depth0.25pt\hss}}
\def\Koordinatenachsen(#1,#2)(#3,#4){\unskip
 \hbox to0pt{\hskip-.5pt \ycoord-#4 \advance\ycoord by1
    \vrule height#2 \Einheit width.5pt depth\ycoord \Einheit}%
 \hbox to0pt{\hskip-1 \Einheit \hskip#3\Einheit 
    \xcoord#1 \advance\xcoord by1 \advance\xcoord by-#3 
    \vrule height0.25pt width\xcoord \Einheit depth0.25pt\hss}}
\def\Gitter(#1,#2){\unskip \xcoord0 \ycoord0 \leavevmode
  \LOOP\ifnum\ycoord<#2
    \loop\ifnum\xcoord<#1
      \raise\ycoord \Einheit\hbox to0pt{\hskip\xcoord 
\Einheit\Punkt\hss}%
      \advance\xcoord by1
    \repeat
    \xcoord0
    \advance\ycoord by1
  \REPEAT}
\def\Gitter(#1,#2)(#3,#4){\unskip \xcoord#3 \ycoord#4 \leavevmode
  \LOOP\ifnum\ycoord<#2
    \loop\ifnum\xcoord<#1
      \raise\ycoord \Einheit\hbox to0pt{\hskip\xcoord 
\Einheit\Punkt\hss}%
      \advance\xcoord by1
    \repeat
    \xcoord#3
    \advance\ycoord by1
  \REPEAT}
\def\Label#1#2(#3,#4){\unskip \xdim#3 \Einheit \ydim#4 \Einheit
  \def\lo{\advance\xdim by-.5 \Einheit \advance\ydim by.5 \Einheit}%
  \def\llo{\advance\xdim by-.25cm \advance\ydim by.5 \Einheit}%
  \def\loo{\advance\xdim by-.5 \Einheit \advance\ydim by.25cm}%
  \def\o{\advance\ydim by.25cm}%
  \def\ro{\advance\xdim by.5 \Einheit \advance\ydim by.5 \Einheit}%
  \def\rro{\advance\xdim by.25cm \advance\ydim by.5 \Einheit}%
  \def\roo{\advance\xdim by.5 \Einheit \advance\ydim by.25cm}%
  \def\l{\advance\xdim by-.30cm}%
  \def\r{\advance\xdim by.30cm}%
  \def\lu{\advance\xdim by-.5 \Einheit \advance\ydim by-.6 \Einheit}%
  \def\llu{\advance\xdim by-.25cm \advance\ydim by-.6 \Einheit}%
  \def\luu{\advance\xdim by-.5 \Einheit \advance\ydim by-.30cm}%
  \def\u{\advance\ydim by-.30cm}%
  \def\ru{\advance\xdim by.5 \Einheit \advance\ydim by-.6 \Einheit}%
  \def\rru{\advance\xdim by.25cm \advance\ydim by-.6 \Einheit}%
  \def\ruu{\advance\xdim by.5 \Einheit \advance\ydim by-.30cm}%
  #1\raise\ydim\hbox to0pt{\hskip\xdim
     \vbox to0pt{\vss\hbox to0pt{\hss$#2$\hss}\vss}\hss}%
}
\catcode`\@=12

\usepackage{amsmath,amsthm}
\usepackage{pstricks,pst-node,pstcol}
\usepackage{color}
\usepackage{xspace}
\usepackage[colorlinks=true,
linkcolor=green,
filecolor=brown,
citecolor=green]{hyperref}
\definecolor{gray}{rgb}{.221,.221,.221}

\def\yellow{\textcolor{yellow} }

\def\ep{\epsilon}
\def\is{\ensuremath{\textrm{SeparateIS}}\xspace}
\def\cis{\ensuremath{\textrm{CombineIS}}\xspace}
\def\st{\ensuremath{\textrm{SeparateST}}\xspace}
\def\cst{\ensuremath{\textrm{CombineST}}\xspace}

\begin{document}
\newtheorem{lemma}{Lemma}
\newtheorem{theorem}{Theorem}
\newtheorem*{prop}{Proposition}
\newtheorem{cor}{Corollary}
\begin{center}
{\Large
On Conjugates for Set Partitions and Integer Compositions                        \\ 
}
\vspace{10mm}
DAVID CALLAN  \\
Department of Statistics  \\
\vspace*{-2mm}
University of Wisconsin-Madison  \\
\vspace*{-2mm}
Medical Science Center \\
\vspace*{-2mm}
1300 University Ave  \\
\vspace*{-2mm}
Madison, WI \ 53706-1532  \\
{\bf callan@stat.wisc.edu}  \\
\vspace{5mm}
October 10, 2005
\end{center}

\vspace{3mm}
\begin{center}
   \textbf{Abstract}
\end{center}
There is a familiar conjugate for integer partitions: transpose the Ferrers 
diagram, and a conjugate for integer compositions: transpose a 
Ferrers-like diagram. Here we propose a conjugate for set partitions and 
show that it interchanges \#\,singletons and \#\,adjacencies. Its 
restriction to noncrossing partitions cropped up in a 1972 paper of 
Kreweras. We also exhibit an analogous 
pair of statistics interchanged by the composition conjugate.

\vspace{10mm}

{\Large \textbf{0 \quad The Conjugate of an Integer Partition}  }
A partition of $n$ is a weakly decreasing list of positive integers,
called its parts, whose sum is $n$.
The Ferrers diagram of a partition $a_{1}\ge a_{2}\ge \ldots 
\ge a_{k}\ge 1$ is the $k$-row left-justified array of dots with $a_{i}$ 
dots in the $i$-th row.  The conjugate, obtained by transposing the Ferrers 
diagram, is a well known 
involution on partitions of $n$ that interchanges the largest part and 
the number of parts.

\vspace{7mm}

{\Large \textbf{1 \quad A Conjugate for Set Partitions}  }
The partitions of an $n$-element set, say    
$[n]=\{1,2,\ldots,n\}$, into nonempty blocks are 
counted by the Bell numbers,
\htmladdnormallink{ 
A000110}{http://www.research.att.com:80/cgi-bin/access.cgi/as/njas/sequences/eisA.cgi?Anum= A000110} in OEIS.
A \emph{singleton} is a block containing just 1 element 
and an \emph{adjacency} is an occurrence of two consecutive elements of $[n]$ 
in the same block. Consecutive is used here in the cyclic sense so 
that $n$ and 1 are also considered to be consecutive (the ordinary 
sense is considered in \cite{yang}).
We say $i$ 
\emph{initiates} an adjacency if $i$ and $i+1$\:mod $n$ are in the same block 
and analogously for terminating an adjacency. The number of $k$-block partitions 
of $[n]$ containing \emph{no} adjacencies
is considered in a recent Monthly problem proposed by 
Donald Knuth \cite{affinity}. Suppressing the number of blocks, it 
turns out that \#\,singletons and \#\,adjacencies have the same 
distribution.
 
\begin{theorem}
    There is a bijection $\phi$ on partitions of $[n]$ that 
    interchanges number of singletons and number of adjacencies.
\end{theorem}
To prove this, we need to consider partitions on arbitrary subsets 
rather than just initial segments of the positive integers. (I thank 
Robin Chapman \cite{chapman} for pointing out that my original proof was incorrect.) 
The notions of adjacency, adjacency 
initiator and adjacency terminator generalize in the obvious way.
We always 
write a partition with elements increasing within each block and blocks 
arranged in increasing order of their first (smallest) elements. 
Thus, for example, with a dash separating blocks, the partition 
$\pi =$
3\:5\:12\:-\:4\:8\:10\:-\:7 has support 
supp($\pi)=\{3,4,5,7,8,10,12\}$ and
adjacencies $(8,10),\,(12,3)$.
In a partition with one-element support, this element is considered to 
be both an adjacency initiator and an adjacency terminator.

Consider the operation \is on positive-integer-support partitions defined by 
$\is(\pi)=\big(\rho,(I,S)\big)$ where $I$ is the set of adjacency 
initiators of $\pi$, $S$ is the set of singleton elements of $\pi$, 
and $\rho$ is the partition obtained by suppressing the elements of 
$I\cup S$ in $\pi$.
Thus, for $\pi =$ 3\:5\:12\:-\:4\:8\:10\:-\:7, we have 
$I=\{8,12\},\ S=\{7\}$ and $\rho=$ 3\:5\:-\:4\:10. Also, for a 
one-block partition $\pi=a_{1}a_{2}\ldots a_{k},\ 
\is(\pi)=\big(\epsilon,\,(\{a_{1},\ldots,a_{k}\},\emptyset)\big)$ 
where $\ep$ denotes the empty partition, unless $k=1$ in which case it 
is $\big(\ep,\,(\{a_{1}\},\{a_{1}\})\big)$.
Clearly, for $\rho$ a partition and 
$A,B$ finite sets of positive integers, the  pair $\big(\rho,(A,B)\big)$
lies in the range of \is  if{f}
(i) $\rho=\ep,\ A=B$ and both are singletons, or (ii)
supp$(\rho),\,A,\,B$ are disjoint and for no successive pair $(a,b)$ in 
supp$(\rho)\cup A \cup B$ is $a \in A$ and $b \in B$ (the successor 
of an adjacency initiator cannot be a singleton).

Analogously, define \st with $S$ the set of singleton elements and $T$ 
the set of adjacency terminators. Note that the condition for 
$\big(\rho,(A,B)\big)$ to lie in range(\st) is precisely the same as 
for it to lie in range(\is). Both are injective and so we may define 
their respective inverses \cis and \cst and we will make use of the crucial 
property that these inverses have identical domains. For example, \cst 
is defined on $\big($3\:10\:-\:4\:7\:-\:12\:$,(\{11\},\{1,2\})\big)$ and 
yields the partition 1\:2\:12\:-\:3\:10\:-\:4\:7\:-\:11.

Now we can define the desired bijection. Given a partition $\pi$ on 
$[n]$, form a sequence $\big(\rho_{1},(I_{1},S_{1})\big),\ 
\big(\rho_{2},(I_{2},S_{2})\big),\ \ldots,\ 
\big(\rho_{k},(I_{k},S_{k})\big)$ by setting
$\big(\rho_{1},(I_{1},S_{1})\big)=\is(\pi),\ $ $ 
\big(\rho_{2},(I_{2},S_{2})\big) =\is(\rho_{1}),\ \ldots,\ 
\big(\rho_{k},(I_{k},S_{k})\big) = \is(\rho_{k-1})$ stopping when 
$\rho_{k}$ has no adjacency initiators and no singletons (as must 
eventually occur). For example, with $n=12$ and $\pi=$ 
1\:-\:2\:-\:3\:11\:12\:-\:4\:7\:10\:-\:5\:9\:\:6\:8, the results are 
laid out in the following table
\[
\begin{array}{|c|c|c|c|}\hline
    j & \rho_{j} & I_{j} & S_{j}  \\ \hline \hline
    1 & 3\:12\:-\:4\:7\:10\:-\:5\:9\:-\:6\:8 & \{11\} & \{1,2\}  \\ 
    \hline
    2 & 3\:-\:4\:7\:10\:-\:5\:9\:-\:6\:8 & \{12\} & \emptyset  \\ \hline
    3 & 4\:7\:10\:-\:5\:9\:-\:6\:8 & \emptyset & \{3\}  \\ \hline
    4 & 4\:7\:-\:5\:9\:-\:6\:8 & \{10\} & \emptyset \\ \hline
\end{array}
\]
and $\rho_{4}$ has no adjacency initiators (hence, no adjacencies) and 
no singletons. Next, form a sequence of partitions 
$\tau_{k},\tau_{k-1},\ldots,\tau_{1},\tau_{0}$ by reversing the 
procedure but using \cst rather than \cis. More precisely, set 
$\tau_{k}=\rho_{k},\ \tau_{k-1}=\cst$ on 
$\big(\tau_{k},(I_{k},S_{k})\big),\ \tau_{k-2}=\cst$ on 
$\big(\tau_{k-1},(I_{k-1},S_{k-1})\big),\ \ldots ,
\tau_{0}=\cst$ on $\big(\tau_{1},(I_{1},S_{1})\big)$. Note that for 
$j=k,k-1,\ldots,1$ in turn, \cst is defined on 
$\big(\tau_{j},(I_{j},S_{j})\big)$ because (i) supp($\tau_{j})=$ 
supp($\rho_{j})$, (ii) \cis is certainly defined on  
$\big(\rho_{j},(I_{j},S_{j})\big)$, and (iii) \cst, \cis have the same 
domain. The example yields
\[
\begin{array}{|c|c|c|c|}\hline
    j & \tau_{j} & I_{j} & S_{j}  \\ \hline \hline
    4 & 4\:7\:-\:5\:9\:-\:6\:8 & \{10\} & \emptyset \\ \hline
    3 & 4\:7\:-\:5\:9\:-\:6\:8\:-\:10 & \emptyset & \{3\}  \\ \hline
    2 & 3\:10\:-\:4\:7\:-\:5\:9\:-\:6\:8 & \{12\} & \emptyset  \\ \hline
    1 & 3\:10\:-\:4\:7\:-\:5\:9\:-\:6\:8\:-\:12 & \{11\} & \{1,2\}  \\  
    \hline \hline
    0 & 1\:2\:12\:-\:3\:10\:-\:4\:7\:-\:5\:9\:-\:6\:8\:-\:11 & & \\ \hline
\end{array}
\]
Now $\phi\,:\,\pi \to \tau_{0}$ is the desired bijection: $\phi$ is 
clearly reversible and sends adjacency initiators to singletons and 
singletons to adjacency terminators, and hence interchanges 
\#\,adjacencies and \#\, singletons. \qed

The \emph{complement} of a partition $\pi$ on $[n]$ is $n+1-\pi$ (elementwise).
Our proposed conjugate is as follows.

\noindent\textbf{Definition}\quad \emph{The \emph{conjugate} of a partition $\pi$ 
on $[n]$ is obtained by applying the map $\phi$ 
followed by complementation.}

Since complementation is an involution 
and sends adjacency initiators to adjacency terminators and vice 
versa, it is not hard to deduce the following result.
\begin{theorem}
Conjugation is an involution on partitions of 
$[n]$ that interchanges \# singletons and \# adjacencies.    
\end{theorem}

\vspace*{7mm}

{\Large \textbf{2 \quad Restriction to Noncrossing partitions }  }

A \emph{noncrossing} partition of $[n]$ is one for which 
no quadruple $a<b<c<d$ in $[n]$ has $a,c$ in one block and $b,d$ in 
another. This implies that if the elements of $[n]$ are arranged in a 
circle and neighboring elements 
within each block are joined by line segments, then no line segments 
cross one another in the resulting polygon diagram (see figure below). The bijection $\phi$ of the 
previous section associates with each partition $\pi$ on $[n]$ a partition
$\rho_{k}$ with no adjacency initiators and no singletons.
Since a noncrossing partition cannot avoid both singletons and 
adjacencies, induction yields the following characterization.
\begin{prop}
    $\pi$ is noncrossing if{f} $\rho_{k}$ is the empty partition.
\end{prop}
Both $\phi$ and conjugate preserve the noncrossing property and so we 
have
\begin{prop}
    Conjugation is an involution on noncrossing partitions of 
$[n]$ that interchanges \# singletons and \# adjacencies.   
\end{prop}
On noncrossing partitions, the conjugate coincides with a graphically 
defined bijection 
first considered by Kreweras \cite{krewerasNC} (see also 
\cite{simionNC}). Given a partition $\pi$ on $[n]$, draw its polygon 
diagram. Insert new 
vertices interleaving the old ones and form the partition with largest 
possible blocks subject to its polygon diagram being disjoint from 
that of $\pi$. Label the new vertices clockwise as shown to get 
$\phi(\pi)$ and counterclockwise as shown to get the conjugate of 
$\pi$. 

\vspace*{2mm}

\begin{pspicture}(-2,-2.3)(2.5,2)
\psset{xunit=2cm,yunit=2cm} 
 \SpecialCoor

\pscircle[linecolor=red](0,0){2}

\rput(2.2; 22.5){\textrm{{\footnotesize 2}}}
\rput(2.2; 67.5){\textrm{{\footnotesize 1}}}
\rput(2.2; 112.5){\textrm{{\footnotesize 8}}}
\rput(2.2; 157.5){\textrm{{\footnotesize 7}}}
\rput(2.2; 202.5){\textrm{{\footnotesize 6}}}
\rput(2.2; 247.5){\textrm{{\footnotesize 5}}}
\rput(2.2; 292.5){\textrm{{\footnotesize 4}}}
\rput(2.2; 337.5){\textrm{{\footnotesize 3}}}
\yellow{
\rput(2;22.5){$\bullet$}
\rput(2;337.5){$\bullet$}
\rput(2;292.5){$\bullet$}
}
\rput(2.9;270){\textrm{{\footnotesize 
$\pi=$ 158\,-\,2\,-\,3\,-\,4\,-\,67}}}
\pspolygon[linecolor=yellow](2; 67.5)(2;112.5)(2;247.5)
\psline[linecolor=yellow](2; 157.5)(2;202.5)  
  
\end{pspicture}
\begin{pspicture}(-2.5,-2.3)(2.5,2)
\psset{xunit=2cm,yunit=2cm} 
\SpecialCoor

\pscircle[linecolor=red](0,0){2}

\rput(2.2; 0){\textrm{{\footnotesize 2}}}
\rput(2.2; 45){\textrm{{\footnotesize 1}}}
\rput(2.3; 90){\textrm{{\footnotesize 8}}}
\rput(2.2; 135){\textrm{{\footnotesize 7}}}
\rput(2.3; 180){\textrm{{\footnotesize 6}}}
\rput(2.2; 225){\textrm{{\footnotesize 5}}}
\rput(2.2; 270){\textrm{{\footnotesize 4}}}
\rput(2.2; 315){\textrm{{\footnotesize 3}}}
\yellow{
\rput(2;22.5){$\bullet$}
\rput(2;337.5){$\bullet$}
\rput(2;292.5){$\bullet$}
}
\rput(2.9;270){\textrm{{\footnotesize 
$\phi(\pi)=$ 1234\,-\,57\,-\,6\,-\,8}}}
\pspolygon[linecolor=yellow](2; 67.5)(2;112.5)(2;247.5)
\psline[linecolor=yellow](2; 157.5)(2;202.5) 

\pspolygon[linecolor=black](2; 45)(2;270)(2;315)(2;0)
\psline[linecolor=black](2; 135)(2;225)
\rput(2;90){$\bullet$}
\rput(2;180){$\bullet$}
\end{pspicture}
\begin{pspicture}(-2.5,-2.3)(2.5,2)
\psset{xunit=2cm,yunit=2cm} 
\SpecialCoor

\pscircle[linecolor=red](0,0){2}

\rput(2.2; 0){\textrm{{\footnotesize 7}}}
\rput(2.2; 45){\textrm{{\footnotesize 8}}}
\rput(2.3; 90){\textrm{{\footnotesize 1}}}
\rput(2.2; 135){\textrm{{\footnotesize 2}}}
\rput(2.3; 180){\textrm{{\footnotesize 3}}}
\rput(2.2; 225){\textrm{{\footnotesize 4}}}
\rput(2.2; 270){\textrm{{\footnotesize 5}}}
\rput(2.2; 315){\textrm{{\footnotesize 6}}}
\yellow{
\rput(2;22.5){$\bullet$}
\rput(2;337.5){$\bullet$}
\rput(2;292.5){$\bullet$}
}
\rput(2.9;270){\textrm{{\footnotesize 
conjugate($\pi)=$ 1\,-\,24\,-\,3\,-\,5678}}}
\pspolygon[linecolor=yellow](2; 67.5)(2;112.5)(2;247.5)
\psline[linecolor=yellow](2; 157.5)(2;202.5) 

\pspolygon[linecolor=black](2; 45)(2;270)(2;315)(2;0)
\psline[linecolor=black](2; 135)(2;225)
\rput(2;90){$\bullet$}
\rput(2;180){$\bullet$}
\end{pspicture}

\vspace*{1mm}
 
\begin{center}
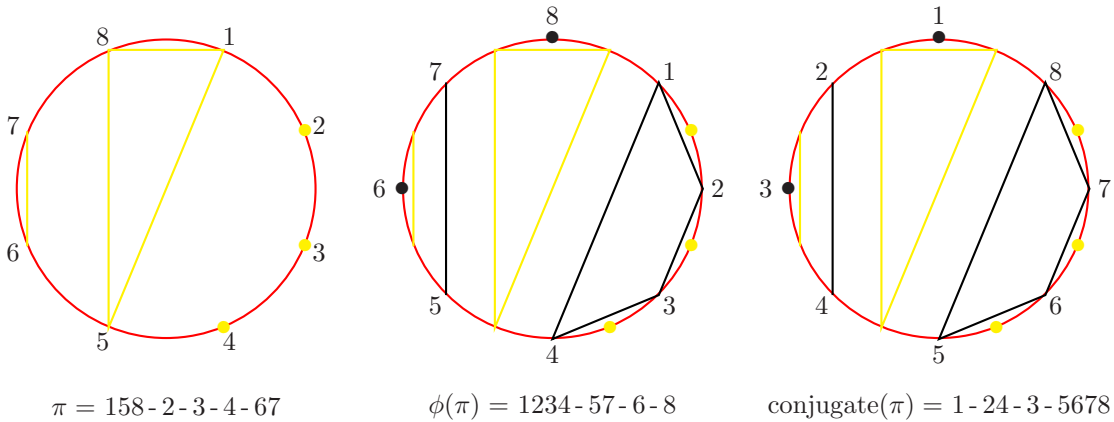
    {\footnotesize polygon diagrams of noncrossing partitions}
\end{center}


{\Large \textbf{3 \quad The Conjugate of an Integer Composition }  }

A \emph{composition} of $n$ is a list 
of positive integers---its parts---whose sum is $n$.
There is a bijection from compositions of $n$ to subsets of $[n-1]$ 
via partial sums:  
$(c_{i})_{i=1}^{k} \mapsto \{ 
\sum_{j=1}^{k}c_{j}\}_{i=1}^{k-1}$, and a further bijection from 
subsets $S$
of $[n-1]$ to lattice paths of $n-1$ unit steps North $(N)$ or East 
($E$): the 
$i$th step is $N$ if $i \in S$ and $E$ otherwise.
The conjugate of a composition is
defined by: pass to lattice path, flip the path in 
the $45^{\circ}$ line, and pass back. For example, with $n=8$ and 
$k=4$,
\[
\overset{\textrm{composition}}{(2,1,2,3)} \rightarrow 
\overset{\textrm{subset}}{\{2,3,5\}}  \rightarrow 
\overset{\textrm{path}}
{ \overset{\textrm{\tiny{1}}}{E}
  \overset{\textrm{\tiny{2}}}{N}
  \overset{\textrm{\tiny{3}}}{N}
  \overset{\textrm{\tiny{4}}}{E}
  \overset{\textrm{\tiny{5}}}{N}
  \overset{\textrm{\tiny{6}}}{E}
  \overset{\textrm{\tiny{7}}}{E} } 
\overset{\raisebox{2.0ex}{\textrm{\scriptsize{f{l}ip}}}}{\ \rightarrow\ }
\overset{\textrm{path}}
{ \overset{\textrm{\tiny{1}}}{N}
  \overset{\textrm{\tiny{2}}}{E}
  \overset{\textrm{\tiny{3}}}{E}
  \overset{\textrm{\tiny{4}}}{N}
  \overset{\textrm{\tiny{5}}}{E}
  \overset{\textrm{\tiny{6}}}{N}
  \overset{\textrm{\tiny{7}}}{N} }  \rightarrow 
\overset{\textrm{subset}}{\{1,4,6,7\}}  \rightarrow  
\overset{\substack{\textrm{conjugate}\\ \textrm{composition}}}{(1,3,2,1,1).}    
\]
There is a neat graphical construction for the lattice path of a 
composition using a kind of shifted Ferrers diagram \cite{combinatoryanalysis}. Represent a part 
$a_{i}$ as a row of $a_{i}$ dots. Stack the rows so each starts where 
its predecessor ends. Then join up the dots with $E$ and $N$ steps.
\vspace*{4mm}
\Einheit=0.6cm
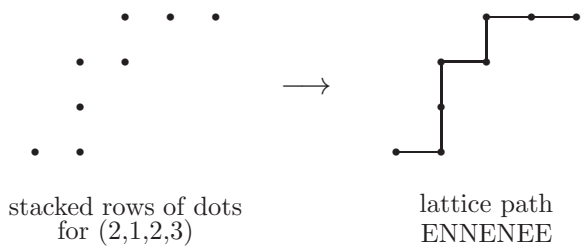
\[
\Pfad(2,2),1221211\endPfad
\DuennPunkt(-6,2)
\DuennPunkt(-5,2)
\DuennPunkt(-5,3)
\DuennPunkt(-5,4)
\DuennPunkt(-4,4)
\DuennPunkt(-4,5)
\DuennPunkt(-3,5)
\DuennPunkt(-2,5)
\DuennPunkt(2,2)
\DuennPunkt(3,2)
\DuennPunkt(3,3)
\DuennPunkt(3,4)
\DuennPunkt(4,4)
\DuennPunkt(4,5)
\DuennPunkt(5,5)
\DuennPunkt(6,5)
\Label\o{\longrightarrow}(0,3)
\Label\o{\textrm{{\footnotesize lattice path}}}(4,0.4)
\Label\o{\textrm{{\footnotesize stacked rows of dots}}}(-4,0.4)
\Label\o{\textrm{{\footnotesize ENNENEE}}}(4,-0.2)
\Label\o{\textrm{{\footnotesize for (2,1,2,3)}}}(-4,-0.2)
\] 
If the compositions of $n$ of a given length are listed in lex 
(dictionary) order, then so are the corresponding subsets, and the 
length of the conjugate is $n+1 -$ length of the original. It follows 
that if the compositions of $n$ are sorted, primarily by 
    length and secondarily by lex order ($n=4$ is shown),
    \[(4),\ (1\  3),\  (2\  2),\  (3\  1),\  (1\  1\  2),\  (1\  2\  1),\  (2\  1\  1),\  (1\  1\  1\  1)
    \]
    then the conjugate of 
    the $i$th composition from the left is the $i$th composition from the 
    right.

There are two statistics on compositions of $n$ 
that are interchanged by conjugation, and these statistics again involve 
``singletons'' (parts $=1$) and ``adjacent'' parts. (For partitions, 
since they
are unordered, adjacency necessarily means ``in value'' but here it 
means ``in position'').
To define them, observe that each part 
has two neighbors except the end parts which have only one or, in the case 
of a one-part composition, none (``wraparound'' neighbors are not allowed 
here). Say a part $\ge 2$ is big; a 
part $=1$ is small. Now define  
\begin{eqnarray*}
    \mu & = &  \textrm{sum of the 
big parts,} \\
   \nu & = & \textrm{sum of the small parts + total number of 
neighbors of the big parts.}
\end{eqnarray*}
For example, in the composition $(3,1,1,4,2)$, 
the big parts 3,\,4,\,2 have 1,\,2,\,1 neighbors respectively; so 
$\mu=3+4+2=9$ and $\nu=(1+1)+(1+2+1)=6$.  
\begin{theorem}
    Conjugation interchanges the statistics $\mu$ and $\nu$ on compositions 
    of $n$ except when $n=1$.
\end{theorem}
\textbf{Proof}\quad Carefully translate $\mu$ and $\nu$ to 
the corresponding lattice paths. 
Using the Iverson 
notation that [\emph{statement}] $=1$ if \emph{statement} is true, 
$=0$ if it is false, we find
\begin{eqnarray*}
    \mu & = & \#\,E\textrm{s}+ \#\,EN\textrm{s} +[\textrm{\,path ends}\ 
    E\,],\  \textrm{and}\\
    \nu & = &  \#\,N\textrm{s}+ \#\,EN\textrm{s} +[\textrm{\,path 
    starts}\  N\,]
\end{eqnarray*}
It is then routine to check that $\mu$ on the flipped path agrees 
with $\nu$ on the original. \qed

Finally, we remark that the genesis of the statistics $\mu$ and $\nu$ 
was the following graphical construction of the conjugate that 
involves ``local'' rather than ``global'' flipping. Suppose given a composition, say $(4,2,1,2,1^{3},3)$, where consecutive 1s have
been collected so that $1^{3}$ is short for $1,1,1$. Represent it
as a list of vertical strips (for parts $>1$) and horizontal
strips (for the 1s):
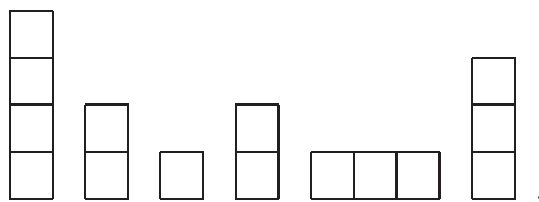
\[
\begin{array}[b]{|c|}
\hline 
\phantom{a}\\ \hline
\phantom{a}\\ \hline
\phantom{a}\\ \hline
\phantom{a}\\ \hline\end{array}
\phantom{aa}
\begin{array}[b]{|c|}
\hline 
\phantom{a}\\ \hline
\phantom{a}\\ \hline\end{array}
\phantom{aa}
\begin{array}[b]{|c|}
\hline 
\phantom{a}\\ \hline\end{array}
\phantom{aa}
\begin{array}[b]{|c|}
\hline 
\phantom{a}\\ \hline
\phantom{a}\\ \hline\end{array}
\phantom{aa}
\begin{array}[b]{|c|c|c|}
\hline 
\phantom{a}&\phantom{a} & \phantom{a} \\ \hline
\end{array}
\phantom{aa}
\begin{array}[b]{|c|}
\hline 
\phantom{a}\\ \hline
\phantom{a}\\ \hline
\phantom{a}\\ \hline\end{array}\ \ .
\]
The vertical strips thus consist of 2 or more squares and 
may occur consecutively, but no two horizontal strips are consecutive.
Insert an (initially) empty horizontal strip between each pair of consecutive 
vertical strips so that horizontal strips $H$ and vertical strips $V$ 
alternate (the first strip may be either an $H$ or a $V$):
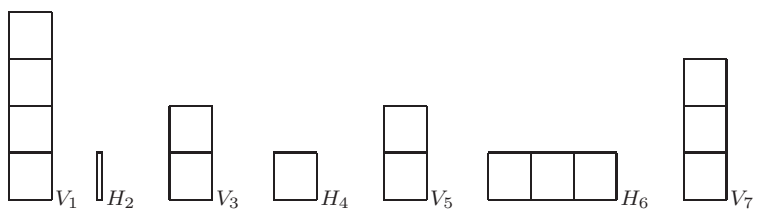
\[
\begin{array}[b]{|c|}
\hline 
\phantom{a}\\ \hline
\phantom{a}\\ \hline
\phantom{a}\\ \hline
\phantom{a}\\ \hline\end{array}_{\,V_{1}} 
\phantom{a}
\begin{array}[b]{|c|}
\hline 
\phantom{\!\!\!\!}\\ \hline
\end{array}_{\,H_{2}} 
\phantom{aa}
\begin{array}[b]{|c|}
\hline 
\phantom{a}\\ \hline
\phantom{a}\\ \hline\end{array}_{\,V_{3}} 
\phantom{aa}
\begin{array}[b]{|c|}
\hline 
\phantom{a}\\ \hline\end{array}_{\,H_{4}} 
\phantom{aa}
\begin{array}[b]{|c|}
\hline 
\phantom{a}\\ \hline
\phantom{a}\\ \hline\end{array}_{\,V_{5}} 
\phantom{aa}
\begin{array}[b]{|c|c|c|}
\hline 
\phantom{a}&\phantom{a} & \phantom{a} \\ \hline
\end{array}_{\,H_{6}} 
\phantom{aa}
\begin{array}[b]{|c|}
\hline 
\phantom{a}\\ \hline
\phantom{a}\\ \hline
\phantom{a}\\ \hline\end{array}_{\,V_{7}} \ \ .
\]
For each horizontal strip $H$, transfer one square from each of its 
neighboring vertical strips to $H$ (there will be two such neighbors in 
general, but possibly just one or even none in the case of the all-1s 
composition):
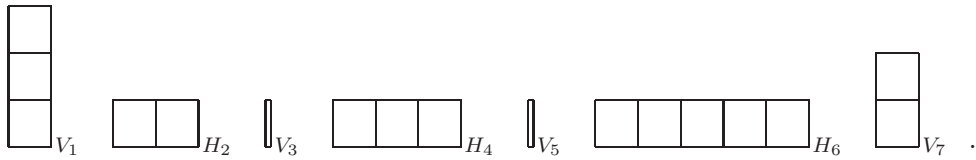
\[
\begin{array}[b]{|c|}
\hline 
\phantom{a}\\ \hline
\phantom{a}\\ \hline
\phantom{a}\\ \hline\end{array}_{\,V_{1}}
\phantom{aa}
\begin{array}[b]{|c|c|}
\hline 
\phantom{a}&\phantom{a}  \\ \hline
\end{array}_{\,H_{2}} 
\phantom{aa}
\begin{array}[b]{|c|}
\hline 
\phantom{\!\!\!\!}\\ \hline
\end{array}_{\,V_{3}}
\phantom{aa}
\begin{array}[b]{|c|c|c|}
\hline 
\phantom{a}&\phantom{a} & \phantom{a} \\ \hline
\end{array}_{\,H_{4}} 
\phantom{aa}
\begin{array}[b]{|c|}
\hline 
\phantom{\!\!\!\!}\\ \hline
\end{array}_{\,V_{5}}
\phantom{aa}
\begin{array}[b]{|c|c|c|c|c|}
\hline 
\phantom{a}&\phantom{a} & \phantom{a}&\phantom{a} & \phantom{a} \\ \hline
\end{array}_{\,H_{6}} 
\phantom{aa}
\begin{array}[b]{|c|}
\hline 
\phantom{a}\\ \hline
\phantom{a}\\ \hline\end{array}_{\,V_{7}}\ \ .
\]
Since each vertical strip originally contained $\ge 2$ squares this 
will always be possible, though some vertical strips may afterward be 
empty, in which case just erase them.  Finally, rotate all strips 
$90^{\circ}$ so that $V$s become $H$s and vice versa:
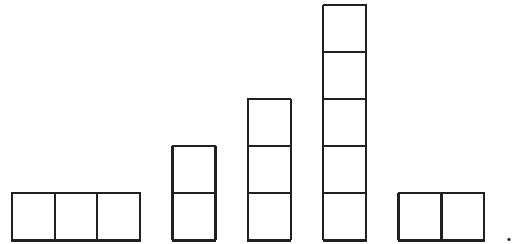
\[
\begin{array}[b]{|c|c|c|}
\hline 
\phantom{a}&\phantom{a} & \phantom{a} \\ \hline
\end{array}
\phantom{aa}
\begin{array}[b]{|c|}
\hline 
\phantom{a}\\ \hline
\phantom{a}\\ \hline\end{array}
\phantom{aa}
\begin{array}[b]{|c|}
\hline 
\phantom{a}\\ \hline
\phantom{a}\\ \hline
\phantom{a}\\ \hline\end{array}
\phantom{aa}
\begin{array}[b]{|c|}
\hline 
\phantom{a}\\ \hline
\phantom{a}\\ \hline
\phantom{a}\\ \hline
\phantom{a}\\ \hline
\phantom{a}\\ \hline\end{array}
\phantom{aa}
\begin{array}[b]{|c|c|}
\hline 
\phantom{a}&\phantom{a}  \\ \hline
\end{array}\ \ .
\phantom{aa}
\]
The result is the conjugate composition: $(1^{3},2,3,5,1^{2})$.

\end{document}